\documentclass{amsart}
\usepackage{amssymb}
\newtheorem{theorem}{Theorem}[section]
\newtheorem{lemma}[theorem]{Lemma}
\newtheorem{proposition}[theorem]{Proposition}

\newtheorem{definition}[theorem]{Definition}
\newtheorem{remark}[theorem]{Remark}
\newcommand{\C}{\mathbb{C}}
\newcommand{\Z}{\mathbb{Z}}
\newcommand{\R}{\mathbb{R}}
\newcommand{\CP}{\mathbb{CP}}
\numberwithin{equation}{section}
\title[Luttinger surgery along Lagrangian tori and non-isotopy]
{Luttinger surgery 
along Lagrangian tori and non-isotopy for singular symplectic plane curves}
\author{D. Auroux}
\address{Department of Mathematics, M.I.T., Cambridge MA 02139, USA}
\email{auroux@math.mit.edu}
\author{S.\,K. Donaldson}
\address{Department of Mathematics, Imperial College, London SW7 2BZ, 
United Kingdom}
\email{s.donaldson@ic.ac.uk}
\author{L. Katzarkov}
\address{Department of Mathematics, University of California, 
Irvine, CA 92697, USA}
\email{lkatzark@math.uci.edu}

\begin{document}
\begin{abstract}
We discuss the properties of a certain type of Dehn surgery along a 
Lagrangian torus in a symplectic 4-manifold, known as Luttinger's surgery,
and use this construction to provide a purely topological interpretation
of a non-isotopy result for symplectic plane curves with cusp and node
singularities due to Moishezon \cite{MChisini}.
\end{abstract}

\maketitle

\section{Introduction}

It is an important open question in symplectic topology to determine
whether, in a given symplectic manifold, all (connected) symplectic 
submanifolds realizing a given homology class are mutually isotopic. 
In the case where the ambient manifold is K\"ahler or complex projective,
one may in particular ask whether symplectic submanifolds are always 
isotopic to complex submanifolds.

The isotopy results known so far rely heavily on the theory of
pseudo-holomorphic curves and on the Gromov compactness theorem
\cite{Gromov}. The best currently known result for smooth curves is due to
Siebert and Tian \cite{ST}, who have proved that smooth connected
symplectic curves of degree at most 17 in $\CP^2$, or realizing homology
classes with intersection pairing at most 7 with the fiber class in a
$S^2$-bundle over $S^2$, are always symplectically isotopic to complex
curves; this strongly suggests that symplectic isotopy holds for smooth
curves in $\CP^2$ and $S^2$-bundles over $S^2$. Isotopy results have also
been obtained for certain singular configurations; e.g., Barraud has
obtained a result for certain arrangements of pseudo-holomorphic
lines in $\CP^2$ \cite{Barraud}.

On the other hand, Fintushel and Stern \cite{FS} and Smith \cite{Smith} 
have constructed infinite families of pairwise non-isotopic smooth connected
symplectic curves representing the same homology classes in certain
symplectic 4-manifolds. In both cases, the construction starts from parallel
copies of a given suitable embedded curve of square zero and modifies them 
by a {\it braiding} construction in order to yield connected symplectic
curves; the constructed submanifolds are distinguished by the diffeomorphism
types of the corresponding double branched covers, either using
Seiberg-Witten theory in the argument of Fintushel and Stern, or by more
topological methods in Smith's argument. It is also worth mentioning that
other examples have recently been obtained by Vidussi using link surgery 
\cite{Vidussi}.

These constructions are predated by a result of Moishezon concerning
singular curves with nodes and cusps in $\CP^2$ \cite{MChisini}. More
precisely, the construction yields infinite families of inequivalent
cuspidal braid monodromies, but as observed by Moishezon the result can be
reformulated in terms of singular plane curves, which one can in fact assume
to be symplectic (cf.\ e.g.\ Theorem 3 of \cite{AK}). The statement can be 
expressed as follows:

\begin{theorem}[Moishezon \cite{MChisini}]
There exists an infinite set $\mathcal{N}$ of positive integers such that,
for each $m\in\mathcal{N}$, there exist integers $\rho_m,d_m$ and
an infinite family of symplectic curves $S_{m,k}\subset\CP^2$
$(k\ge 0)$ of degree $m$ with $\rho_m$ cusps and $d_m$ nodes, such that 
whenever $k_1\neq k_2$ the curves $S_{m,k_1}$ and $S_{m,k_2}$ are not
smoothly isotopic.
\end{theorem}

In particular, because a finiteness result holds for complex curves,
infinitely many of the symplectic curves $S_{m,k}$ are not isotopic to any
complex curve.

Moishezon's argument relies on the observation that the fundamental groups
$\pi_1(\CP^2-S_{m,k})$ are mutually non-isomorphic. However this requires
the heavy machinery of braid monodromy techniques, and in particular the
calculation of the fundamental group of the complement of the branch curve
of a generic polynomial map from $\CP^2$ to itself, carried out in
\cite{MVeronese} (see also \cite{TVeronese}) and preceding papers. The
curious reader is referred to \S 6 of \cite{ADKY} (see also \cite{Tsurv})
for an overview of Moishezon-Teicher braid monodromy techniques.
\medskip

The aim of this paper is to provide a topological interpretation of 
Moishezon's construction, along with an elementary proof of Theorem 1.1;
this reformulation shows that Moishezon's result is very similar to
those of Fintushel-Stern and Smith, in the sense that it also reduces to
a braiding process where the various constructed curves are
distinguished by the topology of associated branched covers. We also
show that these constructions can be thought of in terms of
Luttinger surgery \cite{Luttinger} along Lagrangian tori in a symplectic
4-manifold.

We start by introducing the surgery construction and describing its
elementary properties in \S 2; its interpretation in terms of braiding
constructions for branched covers is discussed in \S 3, while Moishezon's
examples are presented in \S 4.

\section{Luttinger surgery along Lagrangian tori}

\subsection{The surgery construction}

Let $T$ be an embedded Lagrangian torus in a symplectic 4-manifold
$(X,\omega)$, and let $\gamma$ be a simple closed co-oriented loop in $T$.

It is well-known that a neighborhood of $T$ in $X$ can be identified
symplectically with a neighborhood of the zero section in the cotangent
bundle $T^*T\simeq T\times \R^2$ with its standard symplectic structure.
Moreover, $T$ itself can be identified with $\R^2/\Z^2$ in such a way that
$\gamma$ is identified with the first coordinate axis and its
co-orientation coincides with the standard orientation of the
second coordinate axis. Denoting by
$(x_1,x_2)$ the corresponding coordinates on $T$ and by $(y_1,y_2)$ the
dual coordinates in the cotangent fibers, the symplectic form is
given by $\omega=dx_1\wedge dy_1+dx_2\wedge dy_2$. 

Let $r>0$ be such that the set $U_r=(\R^2/\Z^2)\times [-r,r]\times [-r,r]
\subset (\R^2/\Z^2)\times \R^2$ is contained 
in the neighborhood of $T$ over which the identification holds.
Choose a smooth step function $\chi:[-r,r]\to [0,1]$
such that $\chi(t)=0$ for $t\le -\frac{r}{3}$, $\chi(t)=1$ for $t\ge
\frac{r}{3}$, and $\int_{-r}^{\,r} t\,\chi'(t)\,dt=0$ (this last condition
expresses the fact that $\chi$ is {\it centered} around $t=0$).
Given an integer $k\in\Z$, define $\phi_k:U_r-U_{r/2}\to U_r-U_{r/2}$ by
the formulas $\phi_k(x_1,x_2,y_1,y_2)=(x_1+k\chi(y_1),x_2,y_1,y_2)$ if
$y_2\ge \frac{r}{2}$ and $\phi_k=\mathrm{Id}$ otherwise.

Because the support of $\chi$ is contained in $[-\frac{r}{3},\frac{r}{3}]$,
the map $\phi_k$ is actually a diffeomorphism of $U_r-U_{r/2}$; moreover,
$\phi_k$ obviously preserves the symplectic form. Therefore, we can
tentatively make the following definition:

\begin{definition} $X(T,\gamma,k)$ is the manifold obtained from $X$ by
removing a small neighborhood of $T$ and gluing back the standard piece
$U_r$, using the symplectomorphism $\phi_k$ to identify the two sides
near their boundaries. In other terms, $X(T,\gamma,k)=(X-U_{r/2})\cup_{\phi_k} U_r.$
\end{definition}

It can be easily checked that this surgery operation is equivalent to that
introduced by Luttinger in \cite{Luttinger} to study Lagrangian tori in 
$\R^4$ (see also \cite{EP}). 

Forgetting about the symplectic structure, the topological description of
the construction is that of a parametrized $1/k$ Dehn surgery (with
Lagrangian framing): a neighborhood
$T\times D^2$ of $T$ is cut out, and glued back in place by identifying
the two boundaries via a diffeomorphism of $T\times S^1$ that acts trivially
on $H_1(T)$ but maps the homology class of the meridian 
$\mu=\{pt\}\times S^1$ to $[\tilde\mu]=[\mu]+k[\gamma]$.

Observe that the normal bundle to $T$ along $\gamma$ comes equipped with a
natural framing, so that the loop $\gamma$ can be pushed away from $T$ in
a canonical way (up to homotopy), which allows us to define the homotopy
class of $\gamma$ in $\pi_1(X-T)$; comparing the fundamental groups of $X$
and $X(T,\gamma,k)$ with $\pi_1(X-T)$, we see that the surgery
operation preserves the fundamental group (resp.\ first homology group) 
whenever $\gamma^k$ is homotopically (resp.\ homologically) trivial in $X-T$.

The fact that the construction is well-defined symplectically is a
consequence of Moser's stability theorem. More precisely:

\begin{proposition}
$X(T,\gamma,k)$ carries a natural symplectic form $\tilde\omega$,
well-defined up to isotopy independently of the choices made in the 
construction. Moreover, deforming $T$ among Lagrangian tori and
$\gamma\subset T$ by smooth isotopies 
induces a deformation $($pseudo-isotopy$)$ of the symplectic structure 
$\tilde\omega$, and if the symplectic area swept by $\gamma$ is equal to
zero then this deformation preserves the cohomology class $[\tilde\omega]$
and is therefore an isotopy.
\end{proposition}

\proof
Fixing an orientation of $T$ (and therefore of $\gamma$), and
observing that the identification of a neighborhood of $T$ with a
neighborhood of the zero section in $T^*T$ is canonical up to isotopy,
the possible choices for coordinate systems over a neighborhood of $T$ 
differ by isotopies and transformations of the form $(x'_1,x'_2,y'_1,y'_2)=
(x_1+nx_2,x_2,y_1,y_2-ny_1)$ for some integer $n$. 

To handle the case of isotopies, thanks to Moser's stability theorem we only
need to worry about the cohomology class of the symplectic form
$\tilde\omega$ on $\tilde{X}=X(T,\gamma,k)$. Because the surgery affects
only a neighborhood of $T$, once a loop $\delta\subset X-U_r$ homotopic
to the meridian of $T$ in $\tilde{X}$ has been fixed, the cohomology class 
$[\tilde\omega]$ is completely determined by the quantity $\int_D \tilde
\omega$, where $D\subset\tilde{X}$ is a disk such that $\partial D=\delta$ 
and realizing a fixed homotopy class.

Equivalently, if one considers a family depending continuously on a
parameter $t\in [0,1]$, the dependence on $t$ of the cohomology class
$[\tilde\omega_t]$ is exactly given by the symplectic area swept
in $X$ by the meridian loop $\tilde\mu_t=\phi_{k,t}(\partial\Delta)$, where 
$\Delta=\{(0,0,y_1,y_2),\ y_1,y_2\in [-r,r]\}\subset U_r$ (this is because 
$\tilde\mu_t$ bounds a Lagrangian disk $\tilde\Delta_t$ in $\tilde{X}$).
However, observing that the loop $\mu_t=\partial\Delta$, which coincides with
$\tilde\mu_t$ on three of its four sides, bounds a Lagrangian disk
in $X$ and therefore sweeps no area, the symplectic area swept by
$\tilde\mu_t$ is the difference between the area swept by the arc
$\{(0,0,y_1,r),\ y_1\in[-r,r]\}$ and that swept by the arc
$\{(k\chi(y_1),0,y_1,r),\ y_1\in [-r,r]\}$. Using the local expression for
the symplectic form and the fact that the function $\chi$ is centered
($\int_{-r}^r t\,\chi'(t)\,dt=0$), one sees that this is equal to $k$ times
the symplectic area swept by the loop $\{(x_1,0,0,0), x_1\in [0,1]\}$, 
i.e.\ by $\gamma$. Therefore, as long as the loop $\gamma$ is fixed, or
that it is moved in such a way that no symplectic area is swept, we do not
need to worry about continuous deformations of the construction parameters.

Observe that the coordinate change 
$(x'_1,x'_2,y'_1,y'_2)=(x_1+nx_2,x_2,y_1,y_2-ny_1)$ simply amounts to
a modification of the shape of the cut-off region inside each fiber of 
the cotangent bundle (from a square to a parallelogram), which clearly 
has no effect on $\tilde\omega$ (e.g., by the above isotopy argument).
Therefore, to complete the proof we only need to consider the effect of
a simultaneous change of orientation on $T$ and $\gamma$ (recall that
a co-orientation of $\gamma$ inside $T$ is fixed); this simply amounts to
changing $x_1$ and $y_1$ into $-x_1$ and $-y_1$, which clearly does not
affect the construction.
\endproof

Additionally, it is straightforward to check that, if $\gamma^*$ is the loop
$\gamma$ with the opposite co-orientation, then $X(T,\gamma^*,k)$ is
symplectomorphic to $X(T,\gamma,-k)$.
\medskip

\noindent
{\bf Example.} Let $\phi:\Sigma\to\Sigma$ be a symplectomorphism of a 
Riemann surface $(\Sigma,\omega_\Sigma)$, and consider the mapping
torus $Y(\phi)=[0,1]\times \Sigma/(1,x)\sim(0,\phi(x))$. The manifold
$X=S^1\times Y(\phi)$ fibers over $S^1\times S^1$, with monodromy
$\mathrm{Id}$ along the first factor and $\phi$ along the second factor,
and carries a natural symplectic structure $\omega=d\theta\wedge dt + 
\omega_\Sigma$. Let $\gamma$ be any simple closed loop in $\Sigma$. By
picking a point $(\theta_0,t_0)\in S^1\times S^1$, we can embed $\gamma$
as a closed loop $\bar\gamma=\{(\theta_0,t_0)\}\times\gamma$ inside a
fiber of $X$. Observe that $T=S^1\times\{t_0\}\times\gamma$ is an embedded
Lagrangian torus in $(X,\omega)$, containing $\bar{\gamma}$.
It is easy to check that the manifold $X(T,\bar\gamma,k)$ is exactly 
$S^1\times Y(\tau^k\circ \phi)$, where $\tau$ is a Dehn twist about the
loop $\gamma$ (positive or negative depending on the co-orientation).

\subsection{Effect on the canonical class}
We now study the effect of the surgery procedure on the canonical class
$c_1(\tilde{K})$ of $\tilde{X}=X(T,\gamma,k)$. Although there is in general
no natural identification between $H^2(X,\Z)$ and $H^2(\tilde{X},\Z)$ (these
two spaces may even have different ranks), we can compare the two
canonical classes $c_1(K)$ and $c_1(\tilde{K})$ by means of the relative
cohomology groups. Indeed, $H^2(X,T)$ can be identified with $H^2(\tilde{X},T)$
using excision, and we have long exact sequences
$$\cdots\longrightarrow H^1(T)\stackrel{\delta}{\longrightarrow} H^2(X,T)
\stackrel{\iota}{\longrightarrow}
H^2(X)\longrightarrow H^2(T)\longrightarrow\cdots$$
$$\cdots\longrightarrow H^1(T)
\,\smash{\stackrel{\tilde\delta}{\longrightarrow}}\,
H^2(\tilde{X},T) \stackrel{\tilde\iota}{\longrightarrow}
H^2(\tilde{X})\longrightarrow H^2(T)\longrightarrow\cdots$$
The choice of a trivialization $\tau$ of $K$ over $T$
determines a lift $\hat{c}_1(K,\tau)$ of $c_1(K)$
in the relative group $H^2(X,T)$: the relative Chern class
with respect to the chosen trivialization.

Observe that the choice of a section $\sigma$ of the Lagrangian 
Grassmannian $\Lambda(TX)$ over a certain subset of $X$
determines the homotopy class of a trivialization $\tau_\sigma$ of
$K$ over the same subset: indeed, considering a 2-form
$\theta_\sigma$ such that $\mathrm{Ker}\,\theta_\sigma=\sigma$ at every 
point, and given any $\omega$-compatible almost-complex structure $J$, 
the $(2,0)$-component of $\theta$ provides a nowhere vanishing section of 
$K$.
With this understood, we can fix homotopy classes of trivializations $\tau_T$
and $\tilde\tau_T$ of $K$
and $\tilde{K}$ over a neighborhood of $T$ by considering the family of 
Lagrangian tori $(\R^2/\Z^2)\times\{(y_1,y_2)\}$ parallel to $T$ in either
$X$ or $\tilde{X}$ (i.e., the trivialization of the canonical bundle is
given by the $(2,0)$-component of $dy_1\wedge dy_2$). 

Because the trivializations $\tau_T$ and $\tilde\tau_T$ of $K$ and
$\tilde{K}$ coincide over the punctured neighborhood $U_r-U_{r/2}$,
the relative Chern classes $\hat{c}_1(K,\tau_T)$ and $\hat{c}_1(\tilde{K},
\tilde{\tau}_T)$ are equal to each other; therefore, we have
$c_1(\tilde{K})=\tilde\iota(\hat{c}_1(K,\tau_T))$. However, if we consider 
another trivialization $\tau$ of $K_{|T}$, differing from $\tau_T$ by an
element $\nu\in H^1(T)$, then we obtain a different lift 
$\hat{c}_1(K,\tau)=\hat{c}_1(K,\tau_T)+\delta(\nu)$
of $c_1(K)$ in $H^2(X,T)$. It is important to observe that, even though
$\delta(\nu)\in H^2(X,T)$ maps to zero in $H^2(X)$, 
it does not necessarily lie in the kernel of $\tilde\iota:H^2(\tilde{X},T)
\to H^2(\tilde{X})$; in fact, $\tilde\iota(\delta(\nu))$ precisely
measures the obstruction for the trivialization of $\tilde{K}$ determined
by $\tau$ over the subset $U_r-U_{r/2}$ to extend over a neighborhood of 
$T$ in $\tilde{X}$. Therefore, it is easy to check that
$\tilde\iota(\delta(\nu))=-\langle \nu,k[\gamma]\rangle\,PD([T])$, and
hence $c_1(\tilde{K})=\tilde\iota(\hat{c}_1(K,\tau)-\delta(\nu))=
\tilde\iota(\hat{c}_1(K,\tau))+k\langle \nu,[\gamma]\rangle\,PD([T])$.
\medskip

In the special case where there is a proportionality relation of
the form $c_1(K)=\lambda[\omega]$ in $H^2(X,\R)$, it is of particular
interest to study simultaneously the effect of the surgery construction
on the canonical and symplectic classes, by directly considering
$c_1(\tilde{K})-\lambda[\tilde\omega]\in H^2(\tilde{X},\R)$.
The assumption on $c_1(K)$ allows us choose a (Hermitian) connection 
$\nabla$ on the canonical bundle $K$ with 
curvature 2-form $F=-2\pi i\lambda\omega$. Since the surgery only affects
a neighborhood of $T$, we can endow $\tilde{K}$ with a (Hermitian) connection
$\tilde\nabla$ with curvature $\tilde{F}$ that coincides with $\nabla$ 
outside of $\smash{U_{r/2}}$.
Denote as previously by $\tilde\mu=\phi_k(\partial\Delta)$ the meridian of
$T$ in $\tilde{X}$, which bounds a Lagrangian disk $\tilde\Delta$ in 
$\tilde{X}$. Then we have $c_1(\tilde{K})-\lambda[\tilde\omega]=
\alpha\,PD([T])$, where $\alpha=\frac{i}{2\pi}\int_{\tilde\Delta} \tilde{F}$.

We use the above-defined trivializations $\tau_T$ and $\tilde\tau_T$ of $K$
and $\tilde{K}$ over neighborhoods of $T$ in $X$ and $\tilde{X}$. This allows
us to write locally the connection on $K$ in the form
$\nabla=d+2\pi i\lambda(y_1\,dx_1+y_2\,dx_2)+i\beta$, where $\beta$ is a
closed 1-form and can therefore be expressed as $\beta=a_1\,dx_1+
a_2\,dx_2+dh$. The integral of $\tilde{F}$ over $\tilde\Delta$ is given by 
the integral along its boundary $\tilde\mu$ of the 1-form representing the
connection $\tilde\nabla$ in the chosen trivialization of $\tilde{K}$, 
i.e.\ the holonomy of $\tilde\nabla$ along $\tilde\mu$ in the chosen 
trivialization (note that the choice of a homotopy class of trivialization
allows us to view the holonomy as $i\R$-valued rather than $S^1$-valued). 
Since the chosen connections and trivializations of $K$ and 
$\tilde{K}$ coincide along $\tilde\mu$, this is equal to the
holonomy of $\nabla$ along $\tilde\mu$, which is given by the formula
$$i\int_{\tilde\mu} (2\pi\lambda y_1+a_1)\,dx_1+(2\pi\lambda y_2+a_2)\,dx_2+
dh=i\int_{-r}^r (2\pi\lambda y_1+a_1)\,k\chi'(y_1)\,dy_1=ika_1.$$ Observing 
that $ia_1$ is the holonomy of the flat connection $\nabla_{|T}$ along the 
loop $\gamma$ in the given trivialization, we obtain the following:

\begin{definition}
Given a loop $\delta\subset X$ and a homotopy class of trivialization $\tau$
of the canonical bundle $K$ along $\delta$, we define $H(\delta,\tau)$ to be
the real number such that the holonomy of $\nabla$ along $\delta$ is equal
to $-2\pi iH(\delta,\tau)$ in the trivialization $\tau$.
\end{definition}

\begin{proposition}
$c_1(\tilde{K})-\lambda[\tilde\omega]=k\,H(\gamma,\tau_T)\,PD([T])$.
\end{proposition}

\begin{remark} \rm
A change of homotopy class of the trivialization $\tau$ affects
the quantity $H(\delta,\tau)$ by an integer, while a continuous deformation
of the loop $\delta$ affects $H(\delta,\tau)$ by $\lambda$ times the
symplectic area swept.
\end{remark}

\section{Luttinger surgery for branched covers}
In this section, we consider the case where $X$ is a branched cover of
another symplectic 4-manifold $(Y,\omega_Y)$, and show that
the braiding constructions used by Fintushel-Stern \cite{FS} and Smith
\cite{Smith} are a special case of the surgery described in \S 2.

Let $f:X\to Y$ be a covering map with smooth ramification curve $R\subset
X$ and simple branching at the generic points of $R$, such that the branch 
curve $\Sigma=f(R)$ is a symplectic submanifold of
$Y$, immersed except possibly at complex cusp points. Recall that
$X$ carries a natural symplectic structure (up to isotopy), obtained
from the degenerate 2-form $f^*\omega_Y$ by adding a small exact
perturbation along the ramification curve $R$. More precisely, the local
models for $f$ near the points of $R$ allow us to construct an exact 2-form
$\alpha$ such that, at any point of $R$, the restriction of $\alpha$ to
the kernel of $df$ is a positive volume form. The form $\omega=f^*\omega_Y+
\epsilon\alpha$ is then symplectic for $\epsilon>0$ sufficiently small,
and its isotopy class does not depend on the choice of $\alpha$ or
$\epsilon$ (see e.g.\ Proposition 10 of \cite{Au} and Theorem 3
of \cite{AK}).

Consider a loop $a_0$ contained in the smooth part of $\Sigma$, and an
annulus $V_0\subset\Sigma$ forming a neighborhood of $a_0$ in $\Sigma$.
Locally the manifold $Y$ is a fibration over a neighborhood $U$ of the 
origin in $\R^2$, with fibers $V_z\subset Y$ that are smooth symplectic
annuli for all $z\in U$. This fibration carries a natural symplectic 
connection, given by the symplectic orthogonal to $V_z$ at each point.
Given a path $t\mapsto z(t)$ in $U$ starting at the origin, for small
enough values of $t$ we can consider the loops $a_t\subset 
V_{z(t)}$ obtained by parallel transport of $a_0$, and by construction
$\bigcup_{t\ge 0} a_t$ is a smooth Lagrangian
surface in $Y$. 

Assume that, for some value $z_0\in U-\{0\}$, the
symplectic annulus $V_{z_0}$ is contained in the branch curve
$\Sigma$. Assume moreover that, by parallel transport along a certain path
$t\mapsto z(t)$ joining the origin to $z_0$ in $U$, we can construct
an embedded Lagrangian annulus $A=\bigcup_{t\in [0,t_0]} a_t$ in $Y$, 
such that $A\cap\Sigma=\partial A=a_0\cup a_{t_0}\subset V_0\cup V_{z_0}$.
Assume finally that, among the lifts of $A$, exactly two have boundary
contained in the ramification curve $R$ of the map $f$ in $X$; these two 
lifts together form an embedded torus $T\subset X$, and a suitable choice 
of the perturbation $\alpha$ of the pull-back form $f^*\omega_Y$ ensures 
that $T$ is Lagrangian.
Because there is freedom in choosing the local fibration and the path $z(t)$,
the above assumptions can be made to hold in a rather wide range of 
situations, including those considered by Fintushel-Stern and Smith, but
also the examples studied by Moishezon \cite{MChisini}.

Choose a smooth arc $t\mapsto \eta(t)\in a_t$ joining the two boundary
components of $A$, and let $\gamma$ be the loop in $T$ formed by the two
lifts of $\eta$ with end points lying in $R$. Observe that the homotopy class
of the loop $\gamma$ in $T$ does not depend on the choice of the arc $\eta$.
Moreover, an orientation of $a_0$ determines a co-orientation of $\gamma$
in $T$.

We can perform a {\it braiding} construction
on the two parallel annuli $V_0$ and $V_{z_0}$ contained in $\Sigma$,
twisting them $k$ times around each other along the annulus $A$ for any 
given integer $k\in\Z$. The process is described by the following local
model: a neighborhood of $a_0$ in $Y$ is diffeomorphic to $D^2\times 
[-r,r]\times S^1$, where the factor $D^2$ corresponds to $U$ and the factor 
$[-r,r]\times S^1$ corresponds to the annuli $V_z$. The branch curve
$\Sigma$ can be locally identified with the subset $\{\pm u\}\times
[-r,r]\times S^1\subset D^2\times [-r,r]\times S^1$, for some 
$u\in D^2-\{0\}$, while the annulus $A$ corresponds to $[-u,u]\times \{0\}
\times S^1$. Considering a step function $\chi:[-r,r]\to [0,1]$ as in \S 2, 
the twisted curve
$\Sigma(A,k)$ is obtained from $\Sigma$ by replacing $\{\pm u\}\times
[-r,r]\times S^1$ with $\{(\pm u \exp(i\pi k\chi(t)),t),\
t\in[-r,r]\}\times S^1$. Observe that the construction depends on the choice
of an orientation of the factor $[-r,r]$, or equivalently (because $\Sigma$ 
is symplectic) of an orientation of the loop $a_0$.
Moreover, it is easy to check that the construction can
be performed in a way that preserves the symplecticity of the twisted curve.

Recall that the double cover of $D^2$ branched at the two points $\pm u$ is
an annulus $S^1\times [-r,r]$; therefore a local model for $X$ is 
$S^1\times [-r,r]\times [-r,r]\times S^1$, with the torus $T$
corresponding to $S^1\times \{0\}\times\{0\}\times S^1$ and the loop
$\gamma$ corresponding to $S^1\times \{0\}\times \{0\}\times \{pt\}$. 
Also recall that a half-twist exchanging the two points $\pm u$ in $D^2$ lifts 
to a Dehn twist of the annulus $S^1\times [-r,r]$, i.e.\ a transformation
of the form $(x_1,y_1)\mapsto (x_1+\chi(y_1),y_1)$. Therefore, if one 
tries to understand the effect of the twisting construction on the branched 
cover $X$ in terms of cutting out the piece $S^1\times [-r,r]\times
[-r,r]\times S^1$ and gluing it back in place via a nontrivial
diffeomorphism, the difference between the gluing maps on the
two sides $S^1\times [-r,r]\times \{\pm r\}\times S^1$ must be the
$k$-th power of a Dehn twist along the first $S^1$ factor; so if e.g.\ we
take the gluing map to be the identity near $S^1\times [-r,r]\times \{-r\}\times S^1$,
then on the other side it must map $(x_1,y_1,r,x_2)$ to $(x_1+k\chi(y_1),
y_1,r,x_2)$.

The modification undergone by $X$ is therefore exactly
the construction described in \S 2. After checking that the two
constructions also agree from a symplectic point of view, we obtain
the following result:

\begin{proposition}
The branched cover of $Y$ obtained from $X$ by replacing the branch curve 
$\Sigma$ with the twisted curve $\Sigma(A,k)$ is naturally
symplectomorphic to $X(T,\gamma,k)$.
\end{proposition}

Note that the construction only depends on the isotopy classes
of the loop $a_0$ and of the arc $\eta$, even symplectically; this follows from
Proposition 2.2 by observing that, when the arc $\eta$ is deformed in $Y$, the
symplectic area swept by $\gamma$ is always zero (the areas swept
by the two lifts of $\eta$ exactly compensate each other).

\begin{remark} \rm
When the branch curve $\Sigma$ contains $n>2$ parallel annuli $V_{z_i}$, 
one can similarly construct modified symplectic surfaces associated to 
arbitrary elements $b$ of the braid group $B_n$. However, decomposing $b$
into a product of the standard generators of $B_n$ and their inverses 
(or any other half-twists) and starting from a suitable collection of
disjoint Lagrangian annuli in $Y$ with boundary in $\Sigma$, the general 
braiding construction easily reduces to an iteration of the elementary 
process described above. Therefore, assuming that the braid $b$ is
{\em liftable}, i.e.\ compatible with the branching data of the map $f$,
and that it can be decomposed into liftable half-twists
(note that in the case of a double cover all braids are liftable),
we can describe the effect of the general braiding construction on the
symplectic manifold $X$ as a sequence of Luttinger surgeries along disjoint
Lagrangian tori. 
\end{remark}

The examples studied by Fintushel and Stern \cite{FS} show that, in
some cases, the non-triviality of Luttinger surgeries along Lagrangian tori
can be proved using Seiberg-Witten invariants; however in many cases 
it is possible to conclude by much more elementary arguments, as shown 
in \S 4 for Moishezon's examples.
\medskip

We finish this section by observing that, in the context of branched covers,
it is possible to provide a more topological interpretation of the quantity
$H(\gamma,\tau_T)$ introduced in \S 2.2 to describe the effect of 
the surgery on the canonical and symplectic classes of $X$ in the case 
where they are proportional to each other. 

More precisely, assume that
$c_1(K)=\lambda[\omega]$ in $H^2(X,\R)$, where $\omega$ is the symplectic
form induced by a branched covering map $f:X\to Y$ with smooth ramification
curve $R$. Assume moreover that $[\gamma]\in H_1(X,\Z)$ is a torsion
element, i.e.\ $m[\gamma]=0$ for some integer $m\neq 0$, and let $N$ be
a surface with boundary such that $\partial N=\gamma_1\cup\dots\cup
\gamma_m$, where the $\gamma_i$ are parallel copies of $\gamma$ all obtained
as double lifts of arcs in $Y$. Then $f_*N$ is a 2-cycle in $Y$, and we
have the following:

\begin{proposition}
$m\,H(\gamma,\tau_T)=(\lambda[\omega_Y]-c_1(K_Y))\cdot [f_*N]-I(N,R)$, 
where $I(N,R)$ is the algebraic intersection number between $R$ and
$N$, counting the $2m$ intersection points that lie on the boundary of $N$
with multiplicity $1/2$.
\end{proposition}

\proof
By definition, $m\,H(\gamma,\tau_T)=\sum H(\gamma_i,\tau_T)$ can be 
expressed as the difference of two terms, one measuring the integral over $N$
of the curvature of the connection $\nabla$ on the canonical bundle $K_X$,
and the other measuring the obstruction to extending the trivialization of 
$K_X$ given over the boundary of $N$ to a trivialization over all of $N$
(i.e., the {\it relative degree} of $K_X$ over $N$ with respect to the
given boundary trivialization). By assumption, the first term is
proportional to the symplectic area of $N$; observing that the exact
perturbation added to $f^*\omega_Y$ does not contribute to this area (in
fact we could work directly with the degenerate form $f^*\omega_Y$), we
obtain that it is equal to $\lambda[\omega_Y]\cdot [f_*N]$.

In order to compute the relative degree $\deg(K_X,N)$ of $K_X$ over $N$
(the boundary trivialization is implicit in the notation),
we first deform the boundary loops $\gamma_i$ inside $X$ in order to obtain 
loops $\gamma'_i$ bounding a surface $N'$ in $X$, disjoint from $R$, and
$C^1$-close to $\gamma_i$; we can assume that the immersed loops
$f(\gamma'_i)\subset Y$ are in fact embedded. The trivialization $\tau_T$ 
naturally induces a
trivialization of $K_X$ over each loop $\gamma'_i$, and the relative degree
is unaffected by the operation.

Recall that the trivialization $\tau_T$ of the canonical bundle over $\gamma$
is defined by the Lagrangian plane field given by the tangent 
spaces to $T$ along $\gamma$. Deforming to $\gamma'_i$, this
corresponds to a Lagrangian plane field generated by two vector fields,
one tangent to $\gamma'_i$ and the other almost parallel to the normal
direction to $\gamma$ in $T$. This trivialization of $K_X$ is naturally the
lift of a trivialization of $K_Y$ along $f(\gamma'_i)$, determined by two
vector fields, one tangent to $f(\gamma'_i)$ and the other pointing in
a direction transverse to the arc $\eta$ inside the Lagrangian annulus $A$.

Outside of the ramification curve, $K_X$ is isomorphic to $f^*K_Y$, and a
trivialization of $K_Y$ lifts to a trivialization of $K_X$. However, we have
in fact $c_1(K_X)=f^*c_1(K_Y)+PD([R])$, so the relative degree of $K_X$ over
$N'$ differs from that of $K_Y$ over $f_*N'$ by a correction term equal to
the algebraic intersection number of $R$ with $N'$. On the other hand,
the relative degree of $K_Y$ over $f_*N'$ can be evaluated by observing that
each loop $f(\gamma'_i)$ bounds a small disk $D_i$ in $Y$ (because
$f(\gamma'_i)$ is contained in a neighborhood of the arc $\eta$). Moreover,
$f_*N'-\sum D_i$ is a 2-cycle in $Y$, homologous to $f_*N$. Therefore,
$\deg(K_Y,f_*N')=c_1(K_Y)\cdot [f_*N]+\sum \deg(K_Y,D_i)$.

Because of the choice of the trivialization of $K_Y$ along $f(\gamma'_i)$,
it can be checked explicitly that $\deg(K_Y,D_i)$, which measures the
obstruction to extending the trivialization over $D_i$, is equal to $-1$, 
$0$ or $+1$ depending on the chosen perturbation of $\gamma_i$. Moreover,
$\deg(K_Y,D_i)$ is in all cases equal 
to half of the intersection number of $D_i$ with the branch curve $\Sigma$, 
which coincides with the local difference between the intersection numbers
$I(N,R)$ and $I(N',R)$ (once again, counting boundary intersections with
a coefficient $1/2$). Therefore, we have $\sum
\deg(K_Y,D_i)=I(N,R)-I(N',R)$, and so $\deg(K_X,N)=\deg(K_Y,f_*N')+
I(N',R)=c_1(K_Y)\cdot [f_*N]+I(N,R)$, which completes the proof.
\endproof

\section{Non-isotopic singular symplectic plane curves}
\subsection{The manifolds $X_{p,0}$}
Given two symplectic manifolds $Y$ and $Z$, both obtained as
branched covers of the same manifold $M$, and assuming that the branch 
curves $D_g$ and $D_h$ of $g:Y\to M$ and $h:Z\to M$ intersect transversely
in $M$, we can construct a new symplectic manifold $X=Y\times_M Z=
\{(y,z)\in Y\times Z,\ g(y)=h(z)\}$. The manifold $X$ is naturally equipped
with two branched covering structures given by the two projections;
considering e.g.\ the projection onto the first factor, $f:X\to Y$, we
obtain a branched covering map which is simply the pull-back of $h$ via
the map $g$. In particular, the fiber of $f$ above a point $y\in Y$ is
naturally identified to the fiber of $h$ above the point $g(y)\in M$,
the degree of $f$ is equal to that of $h$, and its branch curve is 
$D=g^{-1}(D_h)$.

We consider the case where $M=Y=Z=\CP^2$, and $g:Y\to M$ is a generic
map defined by three polynomials of degree $3$, while $h:Z\to M$ is a
generic map defined by three polynomials of degree $p\ge 2$. We define
$X_{p,0}=Y\times_M Z$, and consider the projection to the first factor, 
$f:X_{p,0}\to Y=\CP^2$, which is a branched covering of degree $p^2$.
It is worth noting that $X_{p,0}$ is in fact a complex surface.
Via a suitable transformation in $PGL(3,\C)$, we can assume
that, outside of a given fixed small ball $B$, the branch curve $D_h$ of the
map $h$ lies arbitrarily close to a union of $d=3p(p-1)$ lines passing
through a single point in $\CP^2$ (observe that $\deg D_h=d$).
The cubic map $g$ can be chosen in such a way that $D_g$ does
not intersect the ball $B$. The branch curve $D=g^{-1}(D_h)$ of $f$ can then
be obtained topologically from the union of $d$ smooth cubics
$C_1,\dots,C_d\subset \CP^2$ lying in a generic pencil, by removing a small
neighborhood of each of the $9$ base points where the $C_i$ intersect and
replacing it with a configuration similar to the branch curve of $h$ (or
rather to $D_h\cap B$). 

The manifold $X_{p,0}$ can also be described as follows. Blow up
$Y=\CP^2$ at the 9 intersection points of the cubics $C_i$, in order to
obtain the rational elliptic surface $E(1)$ with twelve singular fibers and
nine exceptional sections of square $-1$. Let $W$ be the $p^2$-fold 
cover of $E(1)$ branched along $d$ smooth fibers $F_1,\dots,F_d$ (the
proper transforms of the cubics $C_1,\dots,C_d$). The branching pattern 
of the projection $q:W\to E(1)$ is prescribed by that of the map $h$, in the 
sense that $W$ is the pullback of the elliptic fibration
$E(1)$ over $\CP^1$  under the base change map $h_{|S}:S\to\ell=\CP^1$,
where $S$ is the smooth plane curve of degree $p$ obtained as the preimage 
by $h$ of a generic line $\ell\subset\CP^2$. In particular, $W$ is the total
space of an elliptic fibration $\pi$ over the curve $S$ 
of genus $(p-1)(p-2)/2$, with $12p^2$ singular fibers and nine 
exceptional sections $E_1,\dots,E_9$ of square $-p^2$.
We can glue a copy of $\CP^2$ to $W$ along each of the exceptional sections
$E_i$, replacing a neighborhood of $E_i$ with the complement of a smooth
degree $p$ curve in $\CP^2$. It is easy to check that the resulting manifold
$W \#_{\cup E_i} 9\CP^2$ is naturally identical to $X_{p,0}$.

Normalize the Fubini-Study K\"ahler form on $Y=\smash{\CP^2}$ so that its
cohomology class is Poincar\'e dual to the homology class $[L]$ of a line;
the natural symplectic structure $\omega$ induced on $X_{p,0}$ by the 
covering map $f$ is then Poincar\'e dual to the homology class 
$[H]=[f^{-1}(L)]$.

\begin{lemma}
The symplectic and canonical classes of $X_{p,0}$ are related by the
identity $c_1(K)=\lambda_p[\omega]$ in $H^2(X_{p,0},\R)$, where $\lambda_p=
(6p-9)/p$.
\end{lemma}

\proof
The ramification curve $R$ of the branched covering $f:X_{p,0}\to Y=\CP^2$ is
the preimage under the projection to the second factor
$e:X_{p,0}=Y\times_M Z \to Z$ of the ramification curve $R_h$ of the degree
$p$ polynomial map $h:Z\to M$. The curve $R_h$ is a smooth curve of degree 
$3p-3$ in $Z=\CP^2$; in particular, denoting by $[\ell]$ the homology class
of a line in $M=\CP^2$, we have $p[R_h]=(3p-3)[h^{-1}(\ell)]$ in $H_2(Z,\Z)$.
Pulling back by $e$, we obtain the equality
$p[R]=(3p-3)[(h\circ e)^{-1}(\ell)]$ in $H_2(X_{p,0},\Z)$.
Since $h\circ e=g\circ f$, and since $[g^{-1}(\ell)]=3[L]$ in $H_2(Y,\Z)$,
we conclude that $p[R]=(9p-9)[f^{-1}(L)]=(9p-9)[H]$.
Since it is a general fact about branched covers of
$\CP^2$ that $[R]=PD(c_1(K))+3[H]$, we conclude that $p(c_1(K)+3[\omega])=
(9p-9)[\omega]$, or equivalently $c_1(K)=\lambda_p [\omega]$.
\endproof

It is worth noting that, because the symplectic structure on $E(1)$ depends
on the choice of the volumes of the blow-up operations, the symplectic
structure on $W$ depends on the choice of the symplectic areas of the
exceptional sections $E_i$, and is determined only up to deformation
(pseudo-isotopy). The situation that we naturally want to consider is the
limit as the area of $E_i$ and consequently the symplectic volume of the 
copy of $\CP^2$ glued to $W$ along $E_i$ become very small; on the level
of the branch curve $D\subset Y$, this means that the balls around the 
intersection points of the cubics $C_1,\dots,C_d$ that we delete and
replace with copies of $D_h\cap B$ are very small.

\subsection{The manifolds $X_{p,k}$}

The branch curve of $q:W\to E(1)$ consists of $d$ parallel elliptic
curves $F_1,\dots,F_d$ (fibers of $E(1)$), and similarly the branch 
curve of $f:X_{p,0}\to \CP^2$ is obtained from $d$ cubics $C_1,\dots,C_d$
in a pencil by a modification near the base points. Therefore, as discussed
in \S 3 we can construct a Lagrangian annulus $A$ in $E(1)$ (or $\CP^2$) 
that lifts to a Lagrangian torus $T$ in $W$ or $X_{p,0}$, and
Luttinger surgery along $T$ amounts to braiding the branch curve along the
annulus~$A$.

We start with the observation that the $d$ branch points $z_1,\dots,z_d$ of
the simple branched cover $h_{|S}:S\to\CP^1$ can be grouped into pairs of
points with matching branching data; this can be done in many ways, and in
fact amounts to the choice of a degeneration of $S$ to a
nodal curve with $p^2$ rational components intersecting in a total of $d/2$
points. In particular, we can assume that one of the components of the
degenerated curve intersects only once with the others; or equivalently,
we can find two branch points of $h_{|S}$, e.g.\ $z_1$ and $z_2$, and an
arc $\eta_0$ joining them in $\CP^1$, such that the union of two lifts of 
$\eta_0$ forms a closed curve $\gamma_0\subset S$ that separates $S$ into
two components, one of genus $0$ consisting of only one sheet of
$h_{|S}$, and the other of genus $(p-1)(p-2)/2$ consisting of the remaining
$p^2-1$ sheets. Equivalently, observing that $W$ can be constructed from 
$p^2$ copies of the elliptic fibration $E(1)$ by repeatedly performing fiber
sums, $\gamma_0$ can also be thought of as a loop in the base $S$ that 
separates one of the copies of $E(1)$ from the others.

Let $a_0$ be an arbitrary simple closed loop in the fiber $F_1$ above $z_1$,
representing a non-zero homology class in $F_1$ and avoiding the 9 points 
where $F_1$ intersects the exceptional sections of $E(1)$. In fact, the
choice made by Moishezon in \cite{MChisini} amounts to choosing a 
degeneration of the pencil of cubics 
containing $C_1,\dots,C_d$ so that each $C_i$ becomes close to a union of
three lines in $\CP^2$, and taking $a_0$ to be one of the vanishing cycles
for the corresponding degeneration of the fiber $F_1$; but other choices for
$a_0$ are equally suitable. As in \S 3, use parallel transport above the
arc $\eta_0$ to construct a Lagrangian annulus $A\subset E(1)$ joining
$a_0$ to a similar loop in the fiber $F_2$ above $z_2$. Note that we equip
$E(1)$ with a symplectic form which coincides with that of $Y$ outside of
a small neighborhood of the exceptional sections; moreover, we can assume
that $z_1$ and $z_2$ are arbitrarily close to each other, so that the
construction is well-defined and the annulus $A$ remains away from the
exceptional sections. In fact, we could also construct $A$ directly as
an embedded Lagrangian annulus joining the cubics $C_1$ and $C_2$ in $Y=\CP^2$.

By construction, the annulus $A$ lifts to an embedded Lagrangian torus $T$ in 
$W-\bigcup E_i \subset X_{p,0}$, and an embedded arc $\eta\subset A$ which
projects to the arc $\eta_0\subset\CP^1$ lifts to an embedded loop 
$\gamma\subset T$ such that $\pi(\gamma)=\gamma_0\subset S$.
Let $D_{p,k}=D(A,k)$ be the singular plane curve obtained from the branch
curve $D$ of $f$ by twisting $k$ times along the annulus $A$, and let 
$X_{p,k}=X_{p,0}(T,\gamma,k)$ be the symplectic manifold obtained from 
$X_{p,0}$ by twisting $k$ times along the loop $\gamma$ in the Lagrangian
torus $T$. By Proposition 3.1, $X_{p,k}$ is naturally a symplectic branched
cover of $Y=\CP^2$, with branch curve $D_{p,k}$.

Although the description that we give here is very different from that given
by Moishezon in \cite{MChisini}, it is an interesting exercise left to the
reader to check that the two constructions are actually identical. In fact,
because the two loops $\gamma_0$ and $a_0$ can be viewed as vanishing cycles
for degenerations of the base and fiber of $\pi$, the operation of partial 
conjugation of the braid monodromy described by Moishezon exactly amounts to
the braiding construction described in \S 3.

\begin{lemma} 
The homology class $[T]\in H_2(X_{p,k},\Z)$ is not a torsion class.
Moreover, if $p\not\equiv 0 \mod 3$ or $k\equiv 0\mod 3$ then $[T]$ is
primitive.
\end{lemma}

\proof
By Poincar\'e duality, $[T]$ is a non-torsion class if and only if we can
find a $2$-cycle that has non-trivial intersection pairing with $T$; this
is possible if and only if the meridian of $T$ represents a torsion class
in $H_1(X_{p,k}-T,\Z)$. 

Recall from \S 2 that the class $[\tilde\mu]$ of a meridian of $T$ in 
$X_{p,k}$ can be expressed as $[\mu]+k[\gamma]$, where $[\mu]$ is the class
of a meridian of $T$ in $X_{p,0}$; since the complements of
$T$ in $X_{p,0}$ and $X_{p,k}$ are diffeomorphic, it is therefore sufficient
to prove that both $[\mu]$ and $[\gamma]$ are torsion classes in
$H_1(X_{p,0}-T,\Z)$.

We first show that $[\mu]$ is trivial in $H_1(X_{p,0}-T,\Z)$. Consider an arc 
$\xi_0$ in $\CP^1$ which joins the image $z_0$ of a singular fiber of $E(1)$ 
to the branch point $z_1$ of $h_{|S}$ and does not intersect $\eta_0$ in
any other point. Starting from the singular point
in the fiber above $z_0$ and using parallel transport along $\xi_0$,
we can construct a (Lagrangian) disk $D\subset E(1)$, lying above $\xi_0$
and with boundary $\delta$ contained in the smooth fiber $F_1$ above $z_1$
($\delta$ is the vanishing cycle associated to the chosen singular fiber
and the arc $\xi_0$). If the point $z_0$ and the arc $\xi_0$ are 
chosen in a suitable way, we can assume that the intersection number of
$\delta$ with $a_0$ in $F_1$ is equal to $1$ (recall that by assumption
$a_0$ represents a primitive class in $H_1(F_1,\Z)$), and that the disk $D$
does not intersect the exceptional sections of $E(1)$. 
The two lifts via $h_{|S}$ of $\xi_0$ which pass through
the ramification point above $z_1$ form a single arc $\xi$ in $S$
that joins two of the critical values of the elliptic fibration $\pi:W\to S$
and intersects the loop $\gamma_0$ transversely in a single point. Similarly,
the disk $D$ lifts to a sphere (of self-intersection $-2$) in $W-\bigcup E_i$;
by construction, the intersection number of this sphere with the torus $T$
is equal to $1$. Removing a complement of the intersection with $T$
from the sphere, we have realized the meridian $\mu$ as a boundary in
$X_{p,0}-T$, and therefore $[\mu]=0$.

We next consider the loop $\gamma$, which we push away from $T$ by moving
it slightly along the fibers of the elliptic fibration $\pi:W\to S$. In
fact, we can keep moving $\gamma$ along the fibers until it lies in a
neighborhood of one of the exceptional sections $E_i$. Recall that
$\pi(\gamma)=\gamma_0$ bounds a disk $\Delta$ in $S$, corresponding to
one of the sheets of $h_{|S}$; however, the monodromy of the fibration
$\pi$ along $\gamma_0$ is not quite trivial, but differs from identity by
a Dehn twist around each of the nine points where the fiber intersects
the exceptional sections $E_i$. In other terms, the normal bundle to $E_i$,
with its natural trivialization over the boundary $\gamma_0$, has degree
$-1$ over the disk $\Delta$. Therefore, there is an obstruction to
collapsing $\gamma$ inside $W-\bigcup E_i$, but $\gamma$ is homologous to
$-\nu$, where $\nu$ is a small meridian loop around $E_i$ in $W$.

In $X_{p,0}$, a neighborhood of $E_i$ is replaced by the complement
$\CP^2-S$ of a smooth plane curve of degree $p$. Considering a generic line
in $\CP^2$ and removing neighborhoods of its $p$ intersection points with 
$S$, we conclude that $-p[\nu]$ is homologically trivial in $X_{p,0}-T$, 
and therefore that $[\gamma]$ is a torsion element in $H_1(X_{p,0}-T,\Z)$,
which completes the proof that $[T]$ is not torsion in $H_2(X_{p,k},\Z)$.

Another way to look at the loop $\gamma$ is to view $W$ as a fiber sum of
$p^2$ copies of $E(1)$, with the loop $\gamma_0$ in the base $S$ separating
one of the $E(1)$'s from the others. Therefore, $W-\bigcup E_i$ contains
a subset $U$ diffeomorphic to the complement of a fiber and of the 9
exceptional sections in the rational elliptic surface $E(1)$; the loop
$\gamma$ then corresponds to the meridian of the removed fiber in $U$ (with
reversed orientation). However, $U$ can be identified with the complement of
a smooth cubic in $\CP^2$, so $\pi_1(U)=\Z/3$, and therefore $3[\gamma]=0$
in $H_1(X_{p,0}-T,\Z)$.

If $p\not\equiv 0\mod 3$, then we conclude that $[\gamma]=0$, and so
$[\tilde\mu]=[\mu]+k[\gamma]=0$, i.e.\ the meridian $\tilde\mu$ is a 
boundary in $X_{p,k}-T$. Therefore we can find a 2-cycle in $X_{p,k}$ 
which intersects $T$ once, i.e.\ $[T]$ is primitive. When $p$ is a 
multiple of $3$, the same argument holds provided that $k$ is also a
multiple of $3$.
\endproof

\subsection{Proof of Theorem 1.1}
Our strategy to prove Theorem 1.1 is to show that the manifolds $X_{p,k}$
are not symplectomorphic to each other by using Proposition 2.4. We start
with a computation of the quantity $H(\gamma,\tau_T)$ introduced
in \S 2.2 in the case of the Lagrangian torus $T$ and the loop $\gamma$
constructed in \S 4.2:

\begin{lemma} In $X_{p,0}$, we have $H(\gamma,\tau_T)=(2p-3)/p$.
\end{lemma}

\proof
We use Proposition 3.3. The ramification curve $R$ of
$f:X_{p,0}\to Y=\CP^2$ is obtained by gluing together the ramification 
curve of $q:W\to E(1)$, which consists of $d=3p(p-1)$ fibers of $\pi$, 
with the ramification curve of a polynomial map of degree $p$, i.e.\ a 
smooth curve of degree $3p-3$, inside each of the nine copies of $\CP^2$
glued to $W$ along the exceptional sections $E_i$. 

Recall from above that the loop $\gamma$ is homotopic inside $W-\bigcup
E_i$ to a small loop $\bar\nu$ that is the reversed meridian to one of the
exceptional sections; this deformation can be performed without crossing
$R$ (except along $\gamma$ itself). Inside $X_{p,0}$, the loop $\bar\nu$ can
also be viewed as the meridian of the smooth degree $d$ curve $S$ removed
from $\CP^2$ prior to gluing with $W$; therefore, taking $p$ copies of
$\bar\nu$ we obtain the (reversed) boundary of a punctured line in
$\CP^2-S$, which intersects the ramification curve in $3p-3$ points.
Therefore, $p$ copies of $\gamma$ bound a surface $N$ in $X_{p,0}$ such that
$I(N,R)=p-(3p-3)$ (recall that the $2p$ boundary intersections only count
with coefficient $1/2$). Moreover, because the image by $f$ of $\CP^2-S$
is contained in a small ball around one of the base points of the pencil 
of cubics on $Y$, one easily checks that the homology class $[f_*N]$ is
trivial. Therefore we have $p\,H(\gamma,\tau_T)=-I(N,R) = 2p-3$, which
gives the result.

Alternately, remember that $W$ is the fiber sum of $p^2$ copies of $E(1)$,
and so $W-\bigcup E_i$ contains a subset $U$ diffeomorphic to the
complement of a fiber and of the 9 exceptional sections in $E(1)$,
corresponding to one sheet of the branched cover $q:W\to E(1)$. Also,
$\gamma$ corresponds to the meridian of the removed fiber in $U$.
Therefore, three copies of $\gamma$ bound a punctured line $\mathcal{N}$ in 
$U$, which does not intersect the ramification curve anywhere except on the 
boundary, so $I(\mathcal{N},R)=3$; moreover, one easily checks that 
$f_*\mathcal{N}$ has degree 1
in $\CP^2$. Recalling that since $Y=\CP^2$ we have $c_1(K_Y)=-3[\omega_Y]$,
we obtain $3\,H(\gamma,\tau_T)=(\lambda_p+3)-3=\lambda_p=(6p-9)/p$, which
again gives the result.
\endproof

\begin{proposition} For a fixed value of $p\not\equiv 0\mod 3$, the
manifolds $X_{p,k}$ $(k\ge 0)$ are pairwise non-symplectomorphic.
The same result remains true for $p\equiv 0\mod 3$ if we restrict ourselves
to values of $k$ that are multiple of $3$.
\end{proposition}

\proof
The manifolds $X_{p,k}$ are distinguished by the periods of the cohomology
class $\alpha_{p,k}=c_1(K_{X_{p,k}})-\lambda_p [\omega_{X_{p,k}}]$ evaluated
on elements of $H_2(X_{p,k},\Z)$. Indeed,
by Proposition 2.4 and Lemma 4.3 we have $\alpha_{p,k}=k(2p-3)/p\,\, PD([T])$,
and by Lemma 4.2 the homology class $[T]\in H_2(X_{p,k},\Z)$ is primitive,
so the evaluation of $\alpha_{p,k}$ on integer homology classes yields
all integral multiples of $k(2p-3)/p$.
\endproof

In fact, the difference between the branched covering maps $f_{p,k}:X_{p,k}
\to Y=\CP^2$ can be seen on a purely topological level, without considering
symplectic structures. Indeed, defining $[L]$ to be the homology class of
a line in $Y$, the cohomology class of the symplectic form on $X_{p,k}$ is
the Poincar\'e dual of $[f_{p,k}^{-1}(L)]$; and the canonical class of
$X_{p,k}$ is related to the homology class of the ramification curve $R$ of
$f_{p,k}$ by the formula $[R_{p,k}]=c_1(K_{X_{p,k}})+3[f_{p,k}^{-1}(L)]$.
Therefore, the cohomology class $\alpha_{p,k}$ is in fact a
smooth invariant of the branched covering structure, and the maps 
$f_{p,k}:X_{p,k}\to\CP^2$ are not even smoothly isotopic as branched covers.

The branch curves $D_{p,k}$ are symplectic curves of degree $m=3d=9p(p-1)$
in $\CP^2$, and by construction they all have the same numbers of nodes
and cusps (in fact there are
$27(p-1)(4p-5)$ cusps and $27(p-1)(p-2)(3p^2+3p-8)/2$ nodes, as can be
checked e.g.\ using the Pl\"ucker formulas, cf.\ \cite{MChisini}).

In order to conclude that the curves $D_{p,k}$ are not smoothly isotopic,
we need to study the possible $p^2$-fold covers of $\CP^2$ branched
along $D_{p,k}$. These are given by homomorphisms from the fundamental
group $\pi_1(\CP^2-D_{p,k})$ to the symmetric group $S_{p^2}$, satisfying
certain compatibility relations. Because $\pi_1(\CP^2-D_{p,k})$ is finitely
generated and $S_{p^2}$ is a finite group, there are only finitely many
such morphisms, i.e.\ $\CP^2$ admits only finitely many $p^2$-fold 
covers branched over $D_{p,k}$. Because we have infinitely many inequivalent
branched covers $X_{p,k}$, we conclude that infinitely many of the curves
$D_{p,k}$ are not smoothly isotopic. This completes the proof of Theorem 1.1.

\begin{remark}\rm
The number of $p^2$-fold covers of $\CP^2$ branched above $D_{p,k}$ can be 
bounded explicitly by observing that $\pi_1(\CP^2-D_{p,k})$ is generated by 
$m=\deg D_{p,k}$ small meridian loops, all of which must be mapped to 
transpositions in $S_{p^2}$. However, the structure of 
$\pi_1(\CP^2-D_{p,k})$, as described by Moishezon in \cite{MChisini} using
braid monodromy techniques, implies that there is in fact only one possible
branched covering structure for each of the curves $D_{p,k}$ as soon as 
$p\ge 3$. It then follows immediately from the non-isotopy of the branched
covers $f_{p,k}:X_{p,k}\to\CP^2$ that the curves $D_{p,k}$ are 
all different.
\end{remark}

\begin{remark}\rm
The fact that the homology class $[T]$ fails to be primitive when 
$p\equiv 0\mod 3$ and $k\not\equiv 0\mod 3$ is directly related to
the first homology groups of the manifolds $X_{p,k}$. Indeed, whereas
it can be easily checked that $H_1(X_{p,k},\Z)$ is trivial whenever $p$ is
not a multiple of $3$, it appears that $H_1(X_{p,0},\Z)=\Z/3$ (generated
e.g.\ by $[\gamma]$ or by $[\nu]$) when $p\equiv 0\mod 3$; as a
consequence, when
$p$ is a multiple of $3$ the group $H_1(X_{p,k},\Z)$ is isomorphic to $\Z/3$
for $k\equiv 0\mod 3$ and trivial otherwise.
\end{remark}

\begin{remark}\rm
The construction presented here can be modified in various manners,
e.g.\ by starting with other pairs of branched covers $g:Y\to M$ and 
$h:Z\to M$, or by twisting the branch curves in different ways. This
potentially leads to many more examples of non-isotopic singular
symplectic curves in symplectic 4-manifolds. However, it remains unknown
whether it is possible to construct examples of non-isotopic
smooth connected symplectic curves representing a homology class of
positive square inside a given compact symplectic 4-manifold.
\end{remark}

\begin{remark}\rm
The relation between our strategy to prove Theorem 1.1 (by comparing the
canonical and symplectic classes of the branched covers $X_{p,k}$) and 
the strategy used by Moishezon in \cite{MChisini} (by comparing the
fundamental groups $\pi_1(\CP^2-D_{p,k})$) becomes more apparent if one
considers the observations and conjectures made in \cite{ADKY} about 
the structure of fundamental groups of branch curve complements.
Indeed, Moishezon's argument relies on a computation showing that, while the
fundamental group $\pi_1(\CP^2-D_{p,0})$ is always infinite, the groups
$\pi_1(\CP^2-D_{p,k})$ are finite as soon as $p\ge 3$ and $k\neq 0$, and
have different ranks for different values of $k$.
On the other hand, Conjecture 1.6 in \cite{ADKY} states that, at least for
``sufficiently ample'' simply connected branched covers of $\CP^2$, the
fundamental group of the complement of the branch curve is directly related
to the numerical properties of the symplectic and canonical classes. In
particular, it follows from Theorem 1.5 in \cite{ADKY} that, if the 
canonical and symplectic classes are proportional to each other, then
the fundamental group of the branch curve complement must be infinite; the
converse implication is conjectured to hold as well (assuming again that
the branched cover is simply connected and ``sufficiently ample'').
The fact that Theorem 1.1 can be proved indifferently by considering
fundamental groups of complements or numerical relations in the homology 
of the branched covers can be considered as additional evidence for these 
conjectures.
\end{remark}

\subsection*{Acknowledgements} The authors wish to thank Ron Stern, Ivan
Smith, Tom Mrowka, Fedor Bogomolov and Miroslav Yotov for helpful comments. 
The first and third authors
wish to thank respectively Ecole Polytechnique and Imperial College for 
their hospitality.

\end{document}